\documentclass[11pt,twoside]{article}

\date{} 

\begin{document} 


\centerline{} 

\centerline{} 

\centerline {\Large{\bf Integrals of products of Hurwitz zeta functions}} 
\vspace{0.2cm}


\centerline{\Large{\bf  via Feynman parametrization and two double}}  
\vspace{0.2cm}


\centerline{\Large{\bf sums of Riemann zeta functions}}  

\centerline{} 

\centerline{}

\centerline{\bf {M.A. Shpot}} 
\vspace{0.25cm}


\centerline{Institute for Condensed Matter Physics, Lviv, Ukraine} 


\centerline{}

\centerline{\bf {R.B. Paris}} 
\vspace{0.25cm}


\centerline{Division of Computing and Mathematics,} 

\centerline{University of Abertay Dundee, Dundee DD1 1HG, UK} 


%
%
\newcommand{\pp}{\vline}
\newcommand{\cc}{\mbox{$c\hspace{-1.2mm}\vline$}}
\newcommand{\CC}{\mbox{$C\hspace{-1.92mm}\vline$}}
\newcommand{\zz}{\mbox{$Z\hspace{-2.2mm}Z$}}
\newcommand{\nn}{\mbox{$I\hspace{-2.2mm}N$}}
\newcommand{\rr}{\mbox{$I\hspace{-2.2mm}R$}}
\newcommand{\pz}{\mbox{$\parindent=0mm$}}
\newcommand{\pe}{\mbox{$\parindent=8mm$}}
\newcommand{\bee}{\begin{equation}}
\newcommand{\ee}{\end{equation}}
\def\f#1#2{\mbox{${\textstyle \frac{#1}{#2}}$}}
\def\dfrac#1#2{\displaystyle{\frac{#1}{#2}}}
\newcommand{\fr}{\frac{1}{2}}
\newcommand{\fs}{\f{1}{2}}
\newcommand{\g}{\Gamma}
\newcommand{\br}{\biggr}
\newcommand{\bl}{\biggl}
\newcommand{\ra}{\rightarrow}
\renewcommand{\topfraction}{0.9}
\renewcommand{\bottomfraction}{0.9}
\renewcommand{\textfraction}{0.05}
\newcommand{\mcol}{\multicolumn}
\begin{abstract}
We consider two integrals over $x\in [0,1]$ involving products of the function $\zeta_1(a,x)\equiv \zeta(a,x)-x^{-a}$, where $\zeta(a,x)$ is the Hurwitz zeta function, given by
$$\int_0^1\zeta_1(a,x)\zeta_1(b,x)\,dx\quad\mbox{and}\quad \int_0^1\zeta_1(a,x)\zeta_1(b,1-x)\,dx$$ when
$\Re (a,b)>1$. These integrals have been investigated recently in \cite{SCP}; here we provide an alternative derivation by application of Feynman parametrization.
We also discuss a moment integral and the evaluation of two doubly infinite sums containing the Riemann zeta function $\zeta(x)$ and two free parameters $a$ and $b$. The limiting forms of these sums when $a+b$ takes on integer values are considered.
\vspace{0.4cm}

\vspace{0.3cm}

\end{abstract}

\vspace{0.3cm}

\noindent{\bf 1.\ Introduction}
\setcounter{section}{1}
\setcounter{equation}{0}
\renewcommand{\theequation}{\arabic{section}.\arabic{equation}}
The Hurwitz zeta function $\zeta(a,x)$
is one of the most fundamental functions in mathematics.
It has important applications, for example, in number theory \cite{Apostol13,Arakawa14},
probability theory \cite{Karatsuba92} and
also in numerous areas of mathematical physics;
see \cite{Elizalde12,KG99,KR04,Krech} for extensive lists of physical examples.
The function $\zeta(a,x)$ is defined by the series
\begin{equation}\label{e11}
\zeta(a,x)=\sum_{k=0}^\infty\frac{1}{(k+x)^a}\qquad(\Re(a)>1;\ x\neq 0, -1, -2, \ldots)
\end{equation}
and elsewhere by analytic continuation,
apart from $a=1$, where it has a simple pole with unit residue.
In its convergence domain, the series (\ref{e11}) converges absolutely and uniformly.
It reduces to the Riemann zeta function when $x=1$, viz.
\[\zeta(a,1)=\zeta(a)=\sum_{k=1}^\infty\frac{1}{k^a}\qquad(\Re(a)>1)\]
and its behaviour is singular as $x\to0^+$ described by $\zeta(a,x)\sim x^{-a}$.

The auxiliary zeta function, obtained by separating the singular zeroth term from $\zeta(a,x)$, is defined by
\bee\label{e12z}
\zeta_1(a,x):=\zeta(a,x)-x^{-a}.
\ee
The advantage of this function is that it is continuous in the interval $[0,1]$. It follows that
\bee\label{e12zs}
\zeta_1(a,x)=\sum_{k=1}^\infty\frac{1}{(k+x)^a}=\zeta(a,x+1)\qquad (\Re (a)>1),
\ee
from which we see that $\zeta_1(a,0)=\zeta(a)$. The integrals involving products of two auxiliary zeta functions given by
\[I(a,b)=\int_0^1 \zeta_1(a,x) \zeta_1(b,x)\,dx\qquad\mbox{and}\qquad
J(a,b)=\int_0^1\zeta_1(a,x)\zeta_1(b,1-x)\,dx\]
for $\Re (a, b)>1$ have been recently investigated in \cite{SCP}.
The evaluation of this type of integral was initiated in \cite{KK} in the study of the mean-value integral $\int_0^1|\zeta_1(a,x)|^2dx$. This problem has been considered further in \cite{JA, KM2, KM, KKY2} in a number-theoretic context when $a=\fs+it$ and $t\ra\pm\infty$.

It was established in \cite{SCP} that, for $\Re (a, b)>1$,
\[I(a,b)=\frac{1}{a+b-1}+\bl\{B(a+b-1,1-a)+B(a+b-1,1-b)\br\}\,\zeta(a+b-1)\]
\begin{equation}\label{eI}
\hspace{2cm}-\sum_{n=0}^\infty\frac{(a)_n}{(1-b)_{n+1}}\{\zeta(a+n)-1\}
-\sum_{n=0}^\infty\frac{(b)_n}{(1-a)_{n+1}}\{\zeta(b+n)-1\}
\end{equation}
and
\begin{eqnarray}\nonumber
J(a,b)&=&B(1-a,1-b)\{\zeta(a+b-1)-1\}
\\\label{eJ}
&&\hspace{0.8cm}-\sum_{n=0}^\infty\frac{(a)_n}{n!}\,\frac{\zeta(a+n)-1}{n+1-b}
-\sum_{n=0}^\infty\frac{(b)_n}{n!}\,\frac{\zeta(b+n)-1}{n+1-a},
\end{eqnarray}
where
\[B(x,y)=\frac{\g(x) \g(y)}{\g(x+y)}\]
is the beta function.
The integral $I(a,b)$ has been given previously in \cite{JA} and {\cite[Corollary 4]{KM} by employing different methods.
The derivation of (\ref{eJ}) was first carried out in $0<\Re (a, b)<1$ and then relied on analytic continuation arguments to establish its validity in $\Re (a, b)>1$.
The left-hand sides of (\ref{eI}) and (\ref{eJ}) present apparent singularities at integer values of $a, b\geq 2$ and require a limiting procedure for these special values; see \cite{SCP} for details.

In the first part of this paper, we show how the integrals $I(a,b)$ and $J(a,b)$ may be evaluated without appeal to analytic continuation arguments by making use of Feynman parametrization.
This consists of employing the integral representations
\bee\label{eFP}
\frac{1}{A^aC^b}=\frac{1}{B(a,b)} \int_0^1\frac{w^{b-1}(1\!-\!w)^{a-1}\,dw}{[A\!+\!w(C\!-\!A)]^{a+b}}
=\frac{2}{B(a,b)} \int_{-1}^1\frac{(1\!-\!w)^{a-1}(1\!+\!w)^{b-1}\,dw}{[A\!+\!C\!+\!w(C\!-\!A)]^{a+b}},
\ee
where the second integral follows from the first by making the change of variable $w\ra (1-w)/2$. The main advantage of Feynman parametrization is that it transforms the product of two (or more)
denominators, as on the left of (\ref{eFP}), into an integral involving a single denominator.
The resulting expression in the integrand contains some linear combination
of the original values $A$ and $C$.
This feature greatly facilitates integrations or summations
over variables or summation indices originally contained in $A$ and $C$,
as we shall explicitly see in Sections 2 and 3. 

This approach was originally invented by R. P. Feynman for efficient
evaluations of ``diagrams" that are named after this scientist.
By now, it is ubiquitous in calculations in quantum field theory
as well as its applications to certain problems of statistical physics.
Corresponding discussions may be found, for example, in such textbooks as
\cite{Amit,Smirnov,ZJ02}. The interested reader may also find in these references
the general expressions of Feynman parametrization applied to any number
of denominators. 
It is hoped that this approach can be successfully exploited to deal with other integrals
or sums of a similar nature.
As possible candidates we note here the mean square integrals reviewed in \cite{Mats00}, integrals involving products of the alternating counterpart to the Hurwitz zeta function \cite{HKK16}
and double Eisenstein series studied very recently in \cite{KN16}.

In the final sections of the paper we discuss the moment integral
\[H_n(a):=\int_0^1x^n \zeta(a,x)\,dx\qquad (n=0, 1, 2, \ldots),\]
where $\Re (a)<n+1$, $a\neq 1$. The particular case of this integral when $a=m$ is a positive integer satisfying $2\leq m\leq n$ is considered. This result is used in the final section to determine the evaluation of two infinite double sums involving the Riemann zeta function and two free parameters.

Before proceeding, we record some necessary preliminary results related to $\zeta(a,x)$ and the Riemann zeta function $\zeta(s)$.
From the Wilton formula (see, for example, \cite[p.~248, Eq.~(7)]{SC} and
references therein) we have the expansion
\bee\label{e12a}
\sum_{k=0}^\infty \frac{(a)_k}{k!}\,\zeta(a+k,b) x^k=\zeta(a,b-x)\qquad (|x|<|b|),
\ee
where $(a)_k=\g(a+k)/\g(a)=a(a+1)\ldots (a+k-1)$ denotes the Pochhammer symbol and, by convention, $(0)_0 = 1$. Substitution of the value $b=1$ in (\ref{e12a}) produces
\bee\label{e12}
\zeta(a,1-x)=\sum_{k=0}^\infty \frac{(a)_k}{k!}\,\zeta(a+k) x^k  \qquad(|x|<1;\ a\neq 1),
\ee
which can be regarded as a Taylor expansion of $\zeta(a,1-x)$ in powers of $x$.
For $n=1, 2, \ldots\,$, we have
the zeta function and its first derivative at negative even integer arguments given by \cite[p.~605]{DLMF}
\bee\label{e14}
\zeta(-2n)=0,\qquad \zeta'(-2n)=\frac{(-)^n (2n)!}{2^{2n+1} \pi^{2n}}\,\zeta(2n+1).
\ee
We shall also make use of the result
\bee\label{e15}
\zeta(1+\epsilon)=\frac{1}{\epsilon}\{1+\gamma_E\epsilon+O(\epsilon^2)\}\qquad (\epsilon\ra 0),
\ee
where $\gamma_E$ is the Euler-Mascheroni constant.

\vspace{0.6cm}

\noindent{\bf 2.\ Feynman parametrization for $I(a,b)$}
\setcounter{section}{2}
\setcounter{equation}{0}
\renewcommand{\theequation}{\arabic{section}.\arabic{equation}}
We consider the integral
\bee\label{e201}
I(a,b)=\int_0^1 \zeta_1(a,x) \zeta_1(b,x)\,dx\qquad (\Re (a, b)>1)
\ee
and let $\gamma:=a+b-1$, where $\Re (\gamma)>1$. From the definition of $\zeta_1(a,x)$ in (\ref{e12zs}) it follows that
\[I(a,b)=\int_0^1 \sum_{k=1}^\infty \sum_{n=1}^\infty \frac{dx}{(k+x)^a(n+x)^b}.\]

In the product of denominators we now apply the Feynman parametrization given by the first expression in (\ref{eFP})
with $A=k+x$ and $C=n+x$ to find
\[\frac{1}{(k+x)^a(n+x)^b}=\frac{1}{B(a,b)} \int_0^1\frac{w^{a-1}(1-w)^{b-1}}{[k+x+w(n-k)]^{\gamma+1}}\,dw.\]
On reversal of the order of integration this gives
\begin{eqnarray*}
I(a,b)&=&\frac{1}{B(a,b)}\int_0^1w^{a-1}(1-w)^{b-1}\bl(\sum_{k=1}^\infty \sum_{n=1}^\infty \int_0^1\frac{dx}{[k+x+w(n-k)]^{\gamma+1}}\br)\,dw\\
&=&\frac{1}{\gamma B(a,b)}\int_0^1w^{a-1}(1-w)^{b-1}\,S(w)\,dw,
\end{eqnarray*}
where, after performing the straightforward integration over $x$,
\begin{eqnarray*}
S(w)&=&\sum_{k=1}^\infty \sum_{n=1}^\infty\bl\{ [k+w(n-k)]^{-\gamma}-[k+1+w(n-k)]^{-\gamma}\br\}\\
&=&\bl\{\sum_{k=1}^\infty \sum_{n=1}^\infty-\sum_{k=2}^\infty \sum_{n=2}^\infty\br\} [k+w(n-k)]^{-\gamma}.
\end{eqnarray*}
A simple manipulation of the above double sum (permissible since $\Re (\gamma)>1$) then yields
\begin{eqnarray*}
S(w)&=&1+\sum_{k=2}^\infty [k+w(1-k)]^{-\gamma}+\sum_{n=2}^\infty [1+w(n-1)]^{-\gamma}\\
&=&1+\sum_{k=1}^\infty (k+1-kw)^{-\gamma}+\sum_{k=1}^\infty (1+kw)^{-\gamma}.
\end{eqnarray*}
Hence we obtain
\[I(a,b)=\frac{1}{\gamma}+\frac{1}{\gamma B(a,b)}\int_0^1 w^{a-1}(1-w)^{b-1}\sum_{k=1}^\infty \bl\{(k+1-kw)^{-\gamma}+(1+kw)^{-\gamma}\br\} dw.\]

If we make the change of integration variable $w\ra 1-w$ in the integral involving the second sum we see that this term corresponds to symmetrization with respect to $a \Leftrightarrow b$. We therefore have
\bee\label{e202}
I(a,b)=\frac{1}{\gamma}+\frac{1}{\gamma B(a,b)} (T_{a,b}+T_{b,a}),
\ee
where
\begin{eqnarray*}
T_{a,b}&:=&\sum_{k=1}^\infty \frac{1}{(k+1)^\gamma} \int_0^1 \frac{w^{a-1}(1-w)^{b-1}}{(1-\frac{kw}{k+1})^\gamma} dw\\
&\,=&B(a,b) \sum_{k=1}^\infty \frac{1}{(k+1)^\gamma}\,{}_2F_1\bl(a,\gamma;\gamma+1;\frac{k}{k+1}\br)
\end{eqnarray*}
by application of \cite[(15.6.1)]{DLMF}, where ${}_2F_1$ denotes the Gauss hypergeometric function. Consequently we obtain from (\ref{e202}) the result
\bee\label{e203}
I(a,b)=\frac{1}{\gamma}+\frac{1}{\gamma}\sum_{k=1}^\infty \frac{1}{(k+1)^\gamma} \bl\{{}_2F_1\bl(a,\gamma;\gamma+1;\frac{k}{k+1}\br)+{}_2F_1\bl(b,\gamma;\gamma+1;\frac{k}{k+1}\br)\br\}.
\ee

The above expression for $I(a,b)$ can be further transformed by use of the fact that \cite[(15.8.4)]{DLMF}
\[\frac{1}{\gamma}\,{}_2F_1\bl(a,\gamma;\gamma+1;\frac{k}{k+1}\br)=B(1-a,\gamma) \bl(\frac{k}{k+1}\br)^{-\gamma}\hspace{3cm}\]
\[\hspace{6cm}-\frac{(k+1)^{a-1}}{1-a}\,{}_2F_1\bl(1,b;2-a;\frac{1}{k+1}\br).\]
Substitution of this last expression into the infinite sum on the right-hand side of (\ref{e203}) then yields
\[\{B(1-a,\gamma)+B(1-b,\gamma)\}\zeta(\gamma)-\frac{1}{1-a}\sum_{k=1}^\infty\frac{1}{(k+1)^b}
\,{}_2F_1\bl(1,b;2-a;\frac{1}{k+1}\br)\hspace{2cm}\]
\[\hspace{6cm} -\frac{1}{1-b}\sum_{k=1}^\infty\frac{1}{(k+1)^a}
\,{}_2F_1\bl(1,a;2-b;\frac{1}{k+1}\br).\]
Expansion of the hypergeometric functions as series in ascending powers of $(k+1)^{-1}$ and evaluation of the sum over $k$ in terms of the zeta function then finally produces the result, when $\Re (a, b)>1$,
\[I(a,b)=\frac{1}{\gamma}+\{B(1-a,\gamma)+B(1-b,\gamma)\}\zeta(\gamma)\]
\bee\label{e204}
-\sum_{n=0}^\infty\frac{(a)_n}{(1-b)_{n+1}}\{\zeta(a+n)-1\}
-\sum_{n=0}^\infty\frac{(b)_n}{(1-a)_{n+1}}\{\zeta(b+n)-1\}
\ee
in agreement with (\ref{eI}).

We briefly contrast the two forms of expansion for $I(a,b)$ in (\ref{e203}) and (\ref{e204}). The series in (\ref{e203}) can be employed directly for integer values of $a$ and $b$, whereas the series in (\ref{e204}) require a separate limiting treatment for these cases; see \cite{SCP}. However, the price to pay for this simplification lies in the rates of convergence of the infinite series involved: for (\ref{e203}) the terms in the series behave like $k^{-\gamma}$ as $k\ra\infty$, whereas in (\ref{e204}) we have the faster decay rate given by $n^{\gamma-1} 2^{-n}$.

\vspace{0.6cm}

\noindent{\bf 3. \ Feynman parametrization for $J(a,b)$}
\setcounter{section}{3}
\setcounter{equation}{0}
\renewcommand{\theequation}{\arabic{section}.\arabic{equation}}
We consider the integral
\bee\label{e301}
J(a,b)=\int_0^1 \zeta_1(a,x) \zeta(b,1-x)\,dx\qquad (\Re (a,b)>1)
\ee
where again we set $\gamma:=a+b-1$. Following the procedure adopted in Section 2, we have
\[J(a,b)=\int_0^1\sum_{k=1}^\infty\sum_{n=1}^\infty \frac{dx}{(k+x)^a(n+1-x)^b},\]
to which we can apply the Feynman parametrization given by the second expression in (\ref{eFP}) with $A=k+x$ and $C=n+1-x$ yielding
\[\frac{1}{(k+x)^a(n+1-x)^b}=\frac{2}{B(a,b)}\int_{-1}^1\frac{(1-w)^{a-1}(1+w)^{b-1}}{[n+k+1+w(n+1-k-2x)]^{\gamma+1}}\,dw.\]
Then, reversal of the order of integration produces
\begin{eqnarray*}
J(a,b)&=&\frac{2}{B(a,b)}\int_{-1}^1(1-w)^{a-1}(1+w)^{b-1}\\
&&\times\bl(\sum_{k=1}^\infty\sum_{n=1}^\infty \int_0^1\frac{dx}{
[n+k+1+w(n+1-k-2x)]^{\gamma+1}}\br)dw\\
&=&\frac{1}{\gamma B(a,b)}\int_{-1}^1(1-w)^{a-1}(1+w)^{b-1}\,S(w)\,\frac{dw}{w},
\end{eqnarray*}
where, upon evaluation of the inner integral over $x$,
\begin{eqnarray*}
S(w)\!\!&=&\!\!\sum_{k=1}^\infty\sum_{n=1}^\infty\bl\{[n+k+1+w(n+1-k)]^{-\gamma}-[n+k+1+w(n-k-1)]^{-\gamma}\br\}\\
&=&\!\!\bl\{\sum_{k=1}^\infty\sum_{n=2}^\infty-\sum_{k=2}^\infty\sum_{n=1}^\infty\br\} [n+k+w(n-k)]^{-\gamma}\\
&=&\!\! \sum_{k=2}^\infty[k+1-w(k-1)]^{-\gamma}-\sum_{n=2}^\infty [n+1+w(n-1)]^{-\gamma}.
\end{eqnarray*}

Since $S(-w)=-S(w)$, we can write
\begin{eqnarray}
J(a,b)\!\!&=&\!\!\frac{1}{\gamma B(a,b)} \int_0^1\{(1-w)^{a-1}(1+w)^{b-1}+(1-w)^{b-1}(1+w)^{a-1}\}\,S(w)\,\frac{dw}{w}\nonumber\\
&=&\!\!\frac{1}{\gamma B(a,b)}\sum_{k=2}^\infty\frac{1}{(k+1)^\gamma} \int_0^\infty\{u^{a-1}(u+2)^{b-1}+u^{b-1}(u+2)^{a-1}\}\nonumber\\
&&\hspace{3.5cm}\times \bl\{\bl(u+\frac{2}{k+1}\br)^{-\gamma}-\bl(u+\frac{2k}{k+1}\br)^{-\gamma}\br\}\,du,\label{e302}
\end{eqnarray}
where we have made the change of integration variable $w\ra 1/(u+1)$. The integral (\ref{e302}) is an improper integral, since we have the difference of two terms with logarithmic divergence as $u\ra\infty$. However, we can regularize each of the terms in (\ref{e302}) by introducing the small parameter $\epsilon>0$ into the exponent $\gamma$, which makes the integrations over $u$ finite, and subsequently consider the limit $\epsilon\ra 0^+$. Then we write
\bee\label{e303}
J(a,b)=\frac{1}{\gamma B(a,b)} \sum_{k=2}^\infty \frac{L_{a,b}(k)+L_{b,a}(k)}{(k+1)^\gamma},
\ee
where
\[L_{a,b}(k):=\lim_{\epsilon\ra 0} \int_0^\infty u^{a-1}(u+2)^{b-1}\bl\{\bl(u+\frac{2}{k+1}\br)^{-\gamma-\epsilon}-\bl(u+\frac{2k}{k+1}\br)^{-\gamma-\epsilon}\br\}\,du.\]

The integrals appearing in $L_{a,b}(k)$ can be evaluated by making use of the result \cite[(2.2.6.64)]{PBM}
\[\int_0^\infty x^{\alpha-1}(x+y)^{-\rho}(x+z)^{-\lambda} dx=z^{-\lambda} y^{\alpha-\rho} B(\alpha, \rho+\lambda-\alpha)\,{}_2F_1\bl(\alpha, \lambda; \rho+\lambda; 1-\frac{y}{z}\br)
\]
when $0<\Re (\alpha)<\Re (\rho+\lambda)$. Putting $\alpha=a$, $\rho=\gamma+\epsilon$, $\lambda=1-b$, $z=2$ and $y=2/(k+1)$, $2k/(k+1)$, we obtain
\[L_{a,b}(k)=\lim_{\epsilon\ra 0} 2^{-\epsilon} B(a,\epsilon)\bl\{\xi^{1-b-\epsilon} {}_2F_1(a,1-b;a+\epsilon;1-\xi)\hspace{3cm}\]
\bee\label{e303a}
\hspace{4cm}-(1-\xi)^{1-b-\epsilon} {}_2F_1(a,1-b;a+\epsilon;\xi)\br\},
\ee
where $$\xi:=1/(k+1).$$ Application of the limits (\ref{a1}) and (\ref{a2}) derived in Appendix A shows that
\[L_{a,b}(k)=\psi(a)-\psi(1-b)+B(a,b-1) \xi^{1-b} {}_2F_1(a,1-b;2-b;\xi)+h_b(\xi)-h_a(\xi),\]
where $\psi(z)=\g'(z)/\g(z)$ and $h_\mu(z):=\sum_{n=1}^\infty (\gamma)_n z^n/((\mu)_n\,n)$.
A permutation of the indices $a$ and $b$ gives an analogous expression for $L_{b,a}(k)$.
Use of the standard result $\psi(z)-\psi(1-z)=-\pi \cot \pi z$ then produces
\[L_{a,b}(k)+L_{b,a}(k)=-\frac{\pi \sin \pi(a+b)}{\sin \pi a \sin \pi b}
+B(a,b-1) \xi^{1-b} {}_2F_1(a,1-b;2-b;\xi)\]
\[\hspace{3cm}+B(b,a-1) \xi^{1-a} {}_2F_1(b,1-a;2-a;\xi),\]
where the sums $h_a(\xi)$ and $h_b(\xi)$ have cancelled.
Thus
\bee\label{e304}
J(a,b)=M_0+\sum_{k=2}^\infty\bl\{\frac{\xi^{a}}{b-1}\, {}_2F_1(a,1-b;2-b;\xi)+\frac{\xi^{b}}{a-1}\, {}_2F_1(b,1-a;2-a;\xi)\br\},
\ee
where we have put
\bee\label{e305}
M_0=-\frac{\pi}{\gamma B(a,b)}\,\frac{\sin \pi(a+b)}{\sin \pi a \sin \pi b}\, \{\zeta(\gamma)-1-2^{-\gamma}\}=B(1-a,1-b)\,\{\zeta(\gamma)-1-2^{-\gamma}\}.
\ee

Expansion of the hypergeometric functions appearing in (\ref{e304}) as power series in $\xi$ then leads to
\[J(a,b)=M_0-\sum_{n=0}^\infty\frac{(a)_n}{n! (n+1-b)}\{\zeta(a+n)-1-2^{-a-n}\}\hspace{3cm}\]
\[\hspace{3cm}-\sum_{n=0}^\infty\frac{(b)_n}{n! (n+1-a)}\{\zeta(b+n)-1-2^{-b-n}\}\]
\[=M_0'-\sum_{n=0}^\infty\frac{(a)_n}{n!}\,\frac{\zeta(a+n)-1}{n+1-b}-\sum_{n=0}^\infty\frac{(b)_n}{n!}\, \frac{\zeta(b+n)-1}{n+1-a},\]
where
\begin{eqnarray*}
M_0'&=&M_0+2^{-a}\sum_{n=0}^\infty \frac{(a)_n (\fs)^n}{n! (n+1-b)}+2^{-b}\sum_{n=0}^\infty \frac{(b)_n (\fs)^n}{n! (n+1-a)}\\
&=&M_0+2^{-\gamma} B(1-a,1-b)=B(1-a,1-b)\,\{\zeta(\gamma)-1\}
\end{eqnarray*}
by (\ref{b1}) and (\ref{e305}). Hence we finally obtain the result, when $\Re (a, b)>1$,
\begin{eqnarray}
J(a,b)&=&B(1-a,1-b)\{\zeta(a+b-1)-1\}
\nonumber\\
&&\hspace{1cm}-\sum_{n=0}^\infty\frac{(a)_n}{n!}\,\frac{\zeta(a+n)-1}{n+1-b}
-\sum_{n=0}^\infty\frac{(b)_n}{n!}\,\frac{\zeta(b+n)-1}{n+1-a},\label{e306}
\end{eqnarray}
which is the expansion stated in (\ref{eJ}).

An alternative representation for $J(a,b)$ can be obtained from (\ref{e304}) by making use of the well-known transformation
\cite[(15.8.1)]{DLMF}
\[{}_2F_1(a,1-b;2-b;\xi)=(1-\xi)^{-a}\,{}_2F_1\bl(1,a;2-b;\frac{\xi}{\xi-1}\br).\]
Substitution of this result in (\ref{e304}) then yields
\begin{eqnarray}
J(a,b)\!\!\!&=&\!\!\!M_0-\sum_{k=2}^\infty \bl\{\frac{k^{-a}}{1-b}\,{}_2F_1(1,a;2-b;-k^{-1})+\frac{k^{-b}}{1-a}\,{}_2F_1(1,b;2-a;-k^{-1})\br\}\nonumber\\
&=&\!\!\!M_0-\sum_{n=0}^\infty \frac{(-1)^n (a)_n}{(1-b)_{n+1}}\,\{\zeta(a+n)-1\}-\sum_{n=0}^\infty \frac{(-1)^n (b)_n}{(1-a)_{n+1}}\,\{\zeta(b+n)-1\},\nonumber\\
&&\label{e307}
\end{eqnarray}
where $M_0$ is given in (\ref{e305}).
It is seen that the infinite sums in (\ref{e307}) are just the alternating versions of the sums appearing in the expansion for $I(a,b)$ in (\ref{e204}).

We observe that the expansion in (\ref{e306}) has required no analytic continuation in its derivation; in \cite{SCP},
the expansion (\ref{e306}) was first derived for $0<\Re (a, b)<1$ and then analytically continued into $\Re (a, b)>1$. The rates of convergence of the series in (\ref{e306}) and (\ref{e307}) are seen to be comparable with the behaviour $n^{c-1}2^{-n}$ as $n\ra\infty$, where $c=a, b$ for the two series in (\ref{e306}) and $c=\gamma$ for (\ref{e307}).
When $a$ and $b$ assume integer values a limiting procedure is required; the case $a=b=m$ ($m=2, 3, \ldots $) is discussed in \cite[Section 4]{SCP}.
\vspace{0.3cm}

\noindent{\bf Remark 1.}\ \ The expansions for $I(a,b)$ and $J(a,b)$ in (\ref{e204}) and (\ref{e307}) have been derived for $\Re (a, b)>1$. These expansions can be continued analytically into $\Re (a, b)<1$ since both sides
of (\ref{e204}) and (\ref{e307}) are analytic functions of $a$ and $b$, except at $a=1$, $b=1$ where there is a double pole \cite{SCP}. As a special case we consider the parameters with the values $a=\fs+it$, $b=\fs-it$, where $t$ is a real variable, so that $\gamma=0$. A straightforward limiting procedure applied to the first two terms of (\ref{e204}) yields
\[\lim_{\gamma\ra0} \frac{1}{\gamma}+\frac{ \g(1+\gamma)}{\gamma}\bl\{\frac{\g(1-a)}{\g(1-a+\gamma)}+\frac{\g(1-b)}{\g(1-b+\gamma)}\br\}\zeta(\gamma)=\gamma_E-\log\,2\pi+\Re\,\psi(\fs+it),\]
where we have used the fact that $\zeta'(0)=-\fs \log\,2\pi$ and $\psi(1)=-\gamma_E$, where $\gamma_E=0.57721\ldots$ is Euler's constant.
Then we obtain
\[I(\fs+it,\fs-it)=\gamma_E-\log\,2\pi+\Re\,\bl\{\psi(\fs+it)-2 \sum_{n=0}^\infty\frac{\zeta(n+\fs+it)-1}{n+\fs+it}\br\}\]
and, from (\ref{e307}) and \cite[(5.4.4)]{DLMF},
\[J(\fs+it,\fs-it)=-\frac{5\pi}{2\cosh \pi t}-2\Re \sum_{n=0}^\infty(-)^n\frac{\zeta(n+\fs+it)-1}{n+\fs+it}~.\]

\vspace{0.6cm}

\noindent{\bf 4.\ Moment integrals}
\setcounter{section}{4}
\setcounter{equation}{0}
\renewcommand{\theequation}{\arabic{section}.\arabic{equation}}
\newtheorem{theorem}{Theorem}
The moments of $\zeta(a,x)$ over the interval $x\in [0,1]$ are given by
\bee\label{e13}
H_n(a):=\int_0^1 x^n \zeta(a,x)\,dx=\left\{\begin{array}{ll}0& (n=0)\\
\displaystyle {n! \sum_{k=0}^\infty \frac{(a)_k \zeta(a+k)}{(n+k+1)!}} & (n=1, 2, \ldots ) \end{array}\right.
\ee
provided $\Re (a)<n+1$, $a\neq 1$. The case $n=0$ was given by Broughan \cite{Br}. The integrals for positive integer $n$ follow in a straightforward manner by use of (\ref{e12}) with $b=1$ and $z= 1-x$, to find that\footnote{It is worth remarking that the result when $n=0$ also follows from the second expression in (\ref{e13}), since this can be written as
$(a-1)^{-1}\sum_{k=1}^\infty \frac{(a-1)_k}{k!}\,\zeta(a-1+k)$,
which vanishes by (\ref{e31}) when $\Re (a)<1$.}

\begin{eqnarray*}
\int_0^1 x^n\zeta(a,x)\,dx &=&\sum_{k=0}^\infty\frac{(a)_k \zeta(a+k)}{k!}\int_0^1 x^n(1-x)^kdx\nonumber\\
&=&n!\sum_{k=0}^\infty\frac{(a)_k \zeta(a+k)}{(n+k+1)!}\qquad (\Re (a)<n+1).
\end{eqnarray*}

This provides an alternative result to that obtained in
\cite[Theorem 3.7]{EM02}, namely
\begin{equation}\label{e36a}
H_n(a)=\int_0^1 x^n\zeta(a,x)\,dx=-\frac{n!}{\g(a)} \sum_{k=1}^n \frac{\g(a-k) \zeta(a-k)}{(n-k+1)!}\qquad (\Re (a)<n+1),
\end{equation}
which was established by an inductive argument.
By comparing the expressions in (\ref{e13}) and (\ref{e36a}) we arrive at the summation formula
\begin{equation}\label{e36b}
\sum_{k=0}^\infty\frac{\g(a+k) \zeta(a+k)}{(n+k+1)!}=
-\sum_{k=1}^n \frac{\g(a-k)\zeta(a - k)}{(n - k + 1)!}
\end{equation}
for $\Re (a)<n+1$. This expresses an infinite sum in terms of a finite sum.
\vspace{0.2cm}

\noindent{\bf Remark 2.}\ \
The equality (\ref{e36b}) is a consequence of the Wilton formula (\ref{e12}). To see this, we first observe that (\ref{e36b}) can be expressed as
\[\sum_{k=-n}^\infty\frac{\g(a+k) \zeta(a+k)}{(n+k+1)!}=\sum_{k=1}^\infty\frac{\g(a\!-\!n\!-\!1\!+\!k)}{k!}\, \zeta(a\!-\!n\!-\!1\!+\!k)=0\]
by moving  the right-hand side of (\ref{e36b}) to the left and
changing the summation index $k\to k-n-1$ to obtain the second form of the sum.
On the other hand, the special case of the Wilton formula (\ref{e12a}) with $b=2$, $x=1$ yields the result \cite[p.~250, Eq.~(20)]{SC}
\[\sum_{k=1}^\infty \frac{(\alpha)_k}{k!}\,\{\zeta(\alpha+k)-1\}=1.\]
Upon replacing $\alpha$ by $-\alpha$, the sum on the left-hand side  can
be separated into two terms when $\Re (\alpha)>0$ to yield
\bee\label{e31}
\sum_{k=1}^\infty \frac{(-\alpha)_k}{k!}\,\zeta(-\alpha+k)=\sum_{k=0}^\infty \frac{(-\alpha)_k}{k!}=0\qquad (\Re (\alpha)>0).
\ee
Identification of $\alpha$ in (\ref{e31}) with $n+1-a$ then establishes the result (\ref{e36b}).
\vspace{0.3cm}

\noindent{4.1.\ {\it The evaluation of $H_n(a)$ for integer $a$.}
\vspace{0.2cm}

\noindent
We first note that for $a=-m$, $m=0, 1, 2, \ldots\,$ and positive integer $n$, the situation is straightforward since we have from (\ref{e36a}), upon replacement of $\g(-m-k)/\g(-m)$ by $(-)^k m!/(m+k)!$, the result
\bee\label{eH}
H_n(-m)=\int_0^1 x^n \zeta(-m,x)\,dx=m! n! \sum_{k=1}^n \frac{(-)^{k-1} \zeta(-m-k)}{(m+k)! (n-k+1)!}.
\ee
Here the terms corresponding to even $m+k$ vanish on account of the trivial zeros of $\zeta(s)$.
\vspace{0.1cm}

Generically, the moment integrals are more easily evaluated using the finite sum in (\ref{e36a}).
However, when $a=m$, where $m$ is a positive integer satisfying $m\in [2,n]$, a limiting procedure has to be applied to the sum on the right-hand side of (\ref{e36a}),
which is not the case with (\ref{e13}).
We split this sum at $k=m-1$ to rewrite (\ref{e36a}),
with $a=m+\epsilon$ ($\epsilon\ra 0$), as
\begin{equation}\label{e36B}
-\frac{\Gamma(m)}{n!}\int_0^1 x^n\zeta(m,x)\,dx:=S_{m-2}+T_{m-1}+S_m
\end{equation}
\[=\sum_{k=1}^{m-2}\frac{\g(m\!-\!k)\zeta(m\!-\!k)}{(n-k+1)!}+
\frac{\g(1\!+\!\epsilon)\zeta(1\!+\!\epsilon)}{(n-m+2)!}+
\sum_{k=m}^n \frac{\g(m\!-\!k\!+\!\epsilon)\zeta(m\!-\!k\!+\!\epsilon)}{(n-k+1)!}+O(\epsilon).\]
Here $S_{m-2}$ contains the first $m-2$ terms of the original sum in (\ref{e36a}).
With $a=m+\epsilon$, the limit $\epsilon\to 0$ is without difficulty in this summand, as well as in the overall factor
$\Gamma(m+\epsilon)/\Gamma(m)=1+O(\epsilon)$ appearing on the left-hand side of (\ref{e36B}). Consequently
we may include their uninteresting $\epsilon$-dependencies in the order term.

In the $(m-1)$th term, $T_{m-1}$, it is the zeta function $\zeta(1+\epsilon)$ that exhibits a
$1/\epsilon$ pole as $\epsilon\to 0$. We have, by (\ref{e15}) and the fact that $\g(z+\epsilon)=\g(z)\{1+\epsilon \psi(z)+O(\epsilon^2)\}$,
\begin{equation}\label{BN}
T_{m-1}:=\frac{\g(1+\epsilon)\zeta(1+\epsilon)}{(N+1)!}=
\frac1{\epsilon\,(N+1)!}\Big\{1+O(\epsilon^2)\Big\},\qquad N:=n-m+1.
\end{equation}
In the remaining sum $S_m$ in (\ref{e36B}), with summation index $m\leq k\leq n$,
the pole $1/\epsilon$ appears in each
gamma function $\g(m-k+\epsilon)$ as $\epsilon\to 0$ (while each zeta function of the same argument remains
finite). The sum of these pole terms must cancel the $1/\epsilon$ contribution of the term
$T_{m-1}$ given in (\ref{BN}), which we shall now establish.

By changing the summation index $k\to s+m$ in $S_m$, we write
\begin{equation}\label{BM}
S_m:=\sum_{k=m}^n \frac{\g(m-k+\epsilon)\zeta(m-k+\epsilon)}{(n-k+1)!}=
\sum_{s=0}^{N-1} \frac{\g(-s+\epsilon)\zeta(-s+\epsilon)}{(N-s)!},
\end{equation}
with $N$ as defined above.
As $\epsilon\to 0$, the numerator in the last sum is
\begin{equation}\label{BS}
\g(-s+\epsilon)\zeta(-s+\epsilon)=
\frac{(-)^s \zeta(-s)}{\epsilon\,s!}\Big\{1+\epsilon\,\psi(s+1)+\epsilon\,\frac{\zeta'(-s)}{\zeta(-s)}
+O(\epsilon^2)\Big\}.
\end{equation}
Substitution of this result in (\ref{BM}) shows that the singular part of $S_m$, which we denote
as $P[S_m]$, is given by
\begin{equation}\label{BMM}
P[S_m]=\frac1\epsilon\sum_{s=0}^{N-1} \frac{(-)^s\zeta(-s)}{s!(N-s)!}=
\frac1{\epsilon N!}\Big[\zeta(0)+\sum_{s=1}^{N-1}(-)^s\bl(\!\!\begin{array}{c}N\\s\end{array}\!\!\br) \zeta(-s)\Big],
\end{equation}
where the numerical coefficients in the last sum are the binomial coefficients.
Now we take into account that $\zeta(0)=-\fs$ and that for $s\geq 1$ the zeta functions
of negative integer argument $-s$ are expressed in terms of the Bernoulli numbers $B_s$ via
$\zeta(-s)=-B_{s+1}/(s+1)$. This leads us (with $k=s+1$) to
\begin{equation}\label{BNN}
P[S_m]=
-\frac1{2\epsilon N!}\Big[1-2\sum_{k=2}^N (-)^{k}\bl(\!\!\begin{array}{c}N\\k-1\end{array}\!\!\br) \frac{B_k}{k}\Big].
\end{equation}
In the last sum the factor $(-)^k$ can be dropped since $B_{2k+1}=0$ for $k\geq 1$. Then,
making use of the evaluation \cite{R}
\[ \sum_{k=2}^N \bl(\!\!\begin{array}{c}N\\k-1\end{array}\!\!\br) \frac{B_k}{k}=\frac{N-1}{2(N+1)}\qquad (N\geq 2),
\]
we see immediately that $P[S_m]=-\{\epsilon (N+1)!\}^{-1}$ and this exactly cancels the
$O(\epsilon^{-1})$ contribution from (\ref{BN}).

The $O(1)$ contribution to (\ref{e36B}) consists of $S_{m-2}$ and the finite part of the remaining two terms
in this equation. Thus we have, by taking into account (\ref{BN}), (\ref{BM}) and (\ref{BS})
\begin{equation}\label{TB}
-\frac{\Gamma(m)}{n!}\int_0^1 x^n\zeta(a,x)\,dx=
S_{m-2}+
\sum_{k=0}^{N-1} \frac{(-)^k}{k!(N-k)!}
\big\{\psi(k+1)\zeta(-k)+\zeta'(-k)\big\}.
\end{equation}
Splitting the $k=0$ term from the sum in (\ref{TB}), using $\zeta'(0)=-\fs \log\,2\pi$ and $\psi(1)=-\gamma_E$, where $\gamma_E$ is Euler's constant, we obtain
\[-\frac{\Gamma(m)}{n!}\int_0^1 x^n\zeta(a,x)\,dx=S_{m-2}
+\frac{\gamma_E-\log\,2\pi}{2N!}+
\sum_{k=1}^{N-1} \frac{(-)^k}{k!(N-k)!}\]
\[\hspace{5cm}\times\{\psi(k+1)\,\zeta(-k)+\zeta'(-k)\}.
\]
This then leads to the following theorem:
\begin{theorem}$\!\!\!.$\ \ Let $m$, $n$ denote positive integers such that $2\leq m\leq n$ and $N=n-m+1$. Then
\[H_n(m)\!=\!\int_0^1 x^n \zeta(m,x)\,dx\!=\!\frac{n!}{\g(m)}\bl\{\frac{1}{N!}\sum_{k=1}^{N-1}(-)^{k-1}\bl(\!\!\begin{array}{c}N\\k\end{array}\!\!\br)\{\psi(k+1) \zeta(-k)+\zeta'(-k)\}\]
\bee\label{e29}
\hspace{4.5cm}+\frac{\log\,2\pi -\gamma_E}{2N!}-\sum_{k=1}^{m-2}\frac{\g(m-k) \zeta(m-k)}{(n-k+1)!}\br\}.
\ee
\end{theorem}
\bigskip


As examples of Theorem 1, when $m=2$ and $n=3$ we find
\[H_3(2)=\int_0^1 x^3\zeta(2,x)\,dx=6\{\zeta(-1)+\zeta'(-1)\}-\gamma_E+\frac{3}{2}\log\,2\pi\]
\[\hspace{1.1cm}=\frac{3}{2} \log\,2\pi - \gamma_E-6\log\,A,\]
upon using the values $\zeta(-1)=-\f{1}{12}$ and $\zeta'(-1)=\f{1}{12}-\log\,A$,
where $A=1.28242\ldots$ is Glaisher's constant;
see \cite[p.~144]{DLMF}.

Similarly, when $m=3$ and $n=4$ we find
\[H_4(3)=\int_0^1x^4 \zeta(3,x)\,dx=3\log\,2\pi-2\gamma_E-12 \log\,A-\frac{\pi^2}{12}.\]
\vspace{0.6cm}

\noindent{\bf 5.\ Two summation theorems involving products of zeta functions}
\setcounter{section}{5}
\setcounter{equation}{0}
\renewcommand{\theequation}{\arabic{section}.\arabic{equation}}
In this section we present two infinite summation theorems involving products of two Riemann zeta functions.
\begin{theorem}$\!\!\!.$\ \ For $\Re (a)<1$, $\Re (b)<1$, we have the summation
\begin{equation}\label{e21}
S_1(a,b)=\sum_{j=0}^\infty\sum_{k=1}^\infty \frac{(a)_j (b)_k}{(j+k+1)!}\,\zeta(a+j) \zeta(b+k) = B(1-a,1-b)\,\zeta(a+b-1),
\end{equation}
where $B(x,y)=\g(x) \g(y)/\g(x+y)$ is the beta function.
\end{theorem}

\noindent{\it Proof.}\ \ From \cite[Eq.~(3.4)]{EM02}, we have the evaluation of the product of Riemann zeta functions with complementary second arguments over the interval $[0,1]$ given by
\bee\label{e22}
J^*(a,b):=\int_0^1 \zeta(a,x) \zeta(b,1-x)\,dx=B(1-a,1-b)\,\zeta(a+b-1).
\ee
In \cite{EM02} this result was established for $\Re (a)\leq 0$, $\Re (b)\leq 0$. However, since $\zeta(a,x)=O(x^{-a})$ as $x\ra 0$ and $\zeta(b,1-x)=O((1-x)^{-b})$ as $x\ra 1$, both sides of (\ref{e22}) are analytic functions of $a$ and $b$ for $\Re (a)<1$, $\Re (b)<1$. Hence (\ref{e22}) holds by analytic continuation for $\Re (a)<1$, $\Re (b)<1$.

Substitution of the expansion (\ref{e12}) into the integrand in (\ref{e22}), followed by term-by-term integration making use of (\ref{e13}), shows that
\begin{eqnarray*}
J^*(a,b)&=&\sum_{k=1}^\infty \frac{(b)_k \zeta(b+k)}{k!}\int_0^1 x^k \zeta(a,x)\,dx\\
&=&\sum_{j=0}^\infty\sum_{k=1}^\infty \frac{(a)_j (b)_k}{(j+k+1)!}\,\zeta(a+j) \zeta(b+k)
\end{eqnarray*}
whence the result in (\ref{e22}) follows. 
\bigskip

\noindent{\bf Remark 3.}\ \ When $a+b=-1, -3, \ldots$ the sum in (\ref{e21}) vanishes on account of (\ref{e14}).
If $a=-m$, $m=0, 1, 2, \ldots$ the sum over $j$ is finite consisting of $m+2$ terms, since $(-m)_j \zeta(-m+j)=
(-)^m m!$ for $j=m+1$ and $0$ for $j\geq m+2$. Similar considerations apply when $b$ is a non-positive integer.
If both $a$ and $b$ are non-positive integers then the sum in (\ref{e21}) is a finite sum and we find
\bee
S_1(-m,-n)=\sum_{j=0}^{m+1}\sum_{k=1}^{n+1} \frac{(-m)_j (-n)_k}{(j+k+1)!}\,\zeta(-m+j) \zeta(-n+k)\]
\[=\frac{m! n!}{(m+n+1)!}\,\zeta(-m-n-1)\qquad(m,n=0, -1, -2, \ldots).\ee
This sum consequently vanishes whenever the integers $m$ and $n$ are of different parity by (\ref{e14}).
\bigskip

\begin{theorem}$\!\!\!.$\ \ For $\Re (a)<1$, $\Re (b)<1$ and $\Re (a+b)<1$, we have the summation
\[S_2(a,b)=\sum_{j=0}^\infty\sum_{k=0}^\infty \frac{(a)_j (b)_k}{j! k!(j+k+1)}\,\zeta(a+j) \zeta(b+k)\hspace{5cm}\]
\begin{equation}\label{e23}
\hspace{3cm}=
\frac{\pi \g(a+b-1) \zeta(a+b-1)}{\g(a) \g(b)}\left(\frac{\sin \pi a+\sin \pi b}{\sin \pi a\,\sin \pi b}\right).
\end{equation}
\end{theorem}

\noindent{\it Proof.}\ \ From \cite[Eq.~(6.4)]{KM}, we have the evaluation of the product of Riemann zeta functions with equal second arguments over the interval $[0,1]$ given by
\[
I^*(a,b)\!:=\!\!\int_0^1 \zeta(a,x) \zeta(b,x)\,dx\!\!=\!\!\g(a+b-1) \zeta(a+b-1) \left\{\frac{\g(1-a)}{\g(b)}+\frac{\g(1-b)}{\g(a)}\right\}\]
\bee=\!\!\frac{\pi \g(a+b-1) \zeta(a+b-1)}{\g(a) \g(b)}\left(\frac{\sin \pi a+\sin \pi b}{\sin \pi a\,\sin \pi b}\right)\label{e24}
\ee
valid for $\Re (a)<1$, $\Re (b)<1$ and $\Re (a+b)<1$. Because of the symmetry of the integrand under the transformation $x\to 1-x$, we can rewrite the integral as
\begin{eqnarray*}
I^*(a,b)&=&\int_0^1 \zeta(a,1-x) \zeta(b,1-x)\,dx=\sum_{k=0}^\infty \frac{(b)_k \zeta(b+k)}{k!} \int_0^1 x^k \zeta(a,1-x)\,dx\\
&=&\sum_{k=0}^\infty\sum_{j=0}^\infty \frac{(a)_j (b)_k}{j! k!}\,\zeta(a+j)\zeta(b+k)\int_0^1x^{j+k}dx\\
&=&\sum_{j=0}^\infty\sum_{k=0}^\infty \frac{(a)_j (b)_k}{j! k!(j+k+1)}\,\zeta(a+j) \zeta(b+k)
\end{eqnarray*}
upon use of (\ref{e12}). Combination of this result with (\ref{e24}) then concludes the proof. 
\vspace{0.3cm}

\noindent 5.1.\ {\it The case of $S_2(a,b)$ when $a+b$ is a non-positive integer.}
\vspace{0.2cm}

\noindent
When $a+b$ is a non-positive integer a limiting process has to be applied to the right-hand side of (\ref{e23}).
First consider even values of $a+b$. Let $a+b=-2m+\epsilon$, where $m=0,1 , 2, \ldots\,$, and consider $\epsilon\ra 0$. Then use of the reflection formula for the gamma function shows that
\begin{eqnarray*}
S_2(a,b)\!\!&=&\!\!-\lim_{\epsilon\ra 0}\frac{\pi \zeta(-2m-1+\epsilon) \g(1-b)}{\g(a) \g(2m+2-\epsilon) \sin \pi a\,\sin \pi\epsilon} \\
&&\!\!\hspace{4cm}\times\{\cos \pi a\,\sin \pi\epsilon+(1-\cos \pi\epsilon) \sin \pi a\}\\
&=&\!\!-\frac{\g(1\!-\!a) \g(2m\!+\!1\!+\!a)}{(2m+1)!}\,\cos \pi a \ \zeta(-2m\!-\!1) \quad(-2m\!-\!1<\Re (a)<1).
\end{eqnarray*}
For odd values of $a+b$, we find in a similar manner, with $a+b=-2m-1+\epsilon$, that
\begin{eqnarray*}
S_2(a,b)&=&\lim_{\epsilon\ra 0}\frac{\pi\zeta(-2m-2+\epsilon) \g(1-b)}{\g(a) \g(2m+3-\epsilon) \sin \pi\epsilon}\bl\{1+\frac{\sin \pi(a-\epsilon)}{\sin \pi a}\br\}\\
&=&\frac{2 (a)_{2m+2}}{(2m+2)!}\,\zeta'(-2m-2)\qquad (-2m\!-\!2<\Re (a)<1),
\end{eqnarray*}
where we have used the fact that $\zeta(-2m-2+\epsilon)=\epsilon\zeta'(-2m-2)+O(\epsilon^2)$. This leads to the following theorem:
\begin{theorem}$\!\!\!.$\ \ For $a+b=0$ or a negative integer, we have
\[S_2(a,b)=\sum_{j=0}^\infty\sum_{k=0}^\infty \frac{(a)_j (b)_k}{j! k!(j+k+1)}\,\zeta(a+j) \zeta(b+k)\hspace{6cm}\]
\begin{equation}\label{e25}
\hspace{1cm}=\left\{\begin{array}{ll}-\dfrac{\g(1-a)\g(2m+1+a)}{(2m+1)!}\,\cos \pi a\ \zeta(-2m-1) & (a+b=-2m)\\
\dfrac{2(a)_{2m+2}}{(2m+2)!}\,\zeta'(-2m-2) & (a+b=-2m\!-\!1),\end{array}\right.
\ee
where the first equality holds for $-2m-1<\Re (a)<1$ and the second holds for $-2m-2<\Re (a)<1$.
\end{theorem}

Thus, for example, we find the evaluations corresponding to $a+b=0$ and $-1$ given respectively by
\[S_2(a, -a)=\frac{\pi a}{12}\,\cot \pi a\qquad (-1<\Re (a)<1)\]
and
\[S_2(a,-1-a)=a(a+1) \zeta'(-2) =\frac{ab}{4\pi^2}\,\zeta(3)\qquad (-2<\Re (a)<1),\]
where we have used the fact that $\zeta(-1)=-\f{1}{12}$ and employed the second relation in (\ref{e14}).
\bigskip

\noindent{\bf Remark 4.}\ \ When $a$ and $b$ are both negative integers $S_2(a,b)$ reduces to a finite double sum.

\vspace{0.6cm}

\noindent{\bf Appendix A: \ Derivation of the limiting form of the terms in $L_{a,b}(k)$ as $\epsilon\ra 0$ }
\setcounter{section}{1}
\setcounter{equation}{0}
\renewcommand{\theequation}{\Alph{section}.\arabic{equation}}
We determine the limiting forms of the two hypergeometric functions with argument $0<\xi<1$
\[\xi^{1-b-\epsilon} {}_2F_1(a,1-b;a+\epsilon;1-\xi)\qquad \mbox{and}\qquad (1-\xi)^{1-b-\epsilon} {}_2F_1(a,1-b;a+\epsilon;\xi)\]
that appear in $L_{a,b}(k)$ in (\ref{e303a}) as the regularization parameter $\epsilon\ra 0$. From \cite[(15.8.4)]{DLMF}, we have
\[{}_2F_1(a,1-b;a+\epsilon;1-\xi)=\frac{\g(a+\epsilon)\g(1-b-\epsilon)}{\g(a) \g(1-b)}\,\xi^{b-1+\epsilon}{}_2F_1(\epsilon,\gamma+\epsilon;b+\epsilon;\xi)\]
\[\hspace{3.5cm}+\frac{\g(a+\epsilon)\g(b-1+\epsilon)}{\g(\epsilon) \g(\gamma+\epsilon)}\,{}_2F_1(a,1-b;2-b;\xi),\]
where we recall that $\gamma=a+b-1$. Since $(\epsilon)_n=\epsilon (n-1)!+O(\epsilon^2)$ then
\bee\label{a0}
{}_2F_1(\epsilon,\gamma+\epsilon;b+\epsilon;\xi)=1+\sum_{n=1}^\infty \frac{(\epsilon)_n (\gamma+\epsilon)_n}{(b+\epsilon)_n n!} \,\xi^n=1+\epsilon h_b(\xi)+O(\epsilon^2),
\ee
where
\[h_\mu(\xi):=\sum_{n=1}^\infty \frac{(\gamma)_n \xi^n}{(\mu)_n n}\qquad (\xi<1).\]
Hence, upon noting that $\g(\alpha+\epsilon)/\g(\alpha)=1+\epsilon \psi(\alpha)+O(\epsilon^2)$, we have
\[\xi^{1-b-\epsilon} {}_2F_1(a,1-b;a+\epsilon;1-\xi)=1+\epsilon\{\psi(a)-\psi(1-b)\}\hspace{2cm}\]
\bee\label{a1}
\hspace{2cm}+\epsilon \xi^{1-b}B(a,b-1)\,{}_2F_1(a,1-b;2-b;\xi)+\epsilon\, h_b(\xi)+O(\epsilon^2).
\ee

From \cite[(15.8.1)]{DLMF}
\[{}_2F_1(a,1-b;a+\epsilon;\xi)=(1-\xi)^{b-1+\epsilon} {}_2F_1(\epsilon,\gamma+\epsilon;a+\epsilon;\xi),\]
so that upon use of (\ref{a0}) we have
\bee\label{a2}
(1-\xi)^{1-b-\epsilon} {}_2F_1(a,1-b;a+\epsilon;\xi)
=1+\epsilon h_a(\xi)+O(\epsilon^2).
\ee

\vspace{0.6cm}

\noindent{\bf Appendix B: \ A sum of two hypergeometric functions}
\setcounter{section}{2}
\setcounter{equation}{0}
\renewcommand{\theequation}{\Alph{section}.\arabic{equation}}
\newtheorem{lemma}{Lemma}
\begin{lemma}$\!\!\!.$\ \
The sums
\bee\label{b1}
2^{-a}\sum_{n=0}^\infty \frac{(a)_n (\fs)^n}{n! (n+1-b)}+2^{-b}\sum_{n=0}^\infty \frac{(b)_n (\fs)^n}{n! (n+1-a)}=2^{-\gamma} B(1-a,1-b),
\ee
where $\gamma=a+b-1$.
\end{lemma}

\noindent {\it Proof.}\ \ The two sums in (\ref{b1}) can be expressed in terms of Gauss hypergeometric functions of argument $\fs$ as
\[\frac{2^{-a}}{1-b}\,{}_2F_1(a,1-b;2-b;\fs)+\frac{2^{-b}}{1-a}\,{}_2F_1(b,1-a;2-a;\fs).\]
Application of \cite[(15.8.4)]{DLMF} to the first hypergeometric function shows that
\[
\frac{2^{-a}}{1-b}\,{}_2F_1(a,1-b;2-b;\fs)=\frac{2^{-a}\g(1-a)\g(1-b)}{\g(1-\gamma)}\,{}_1F_0(1-b;;\fs)\]
\[\hspace{7cm}-\frac{1}{2(1-a)}\,{}_2F_1(1,1-\gamma;2-a;\fs)\]
\[=2^{-\gamma}B(1-a,1-b)-\frac{1}{2(1-a)}\,{}_2F_1(1,1-\gamma;2-a;\fs),
\]
where we have employed the result ${}_1F_0(\alpha;;\fs)=2^\alpha$. Application of \cite[(15.8.1)]{DLMF} to the second hypergeometric function produces
\[\frac{2^{-b}}{1-a}\,{}_2F_1(b,1-a;2-a;\fs)=\frac{1}{2(1-a)}\,{}_2F_1(1,1-\gamma;2-a;\fs).\]
Combination of these two identities then establishes the lemma.

\vspace{0.6cm}


\end{document}